\documentclass{article}%
\usepackage{amsfonts}
\usepackage{amsmath}
\usepackage{amssymb}
\usepackage{graphicx}%
\newtheorem{theorem}{Theorem}

\newtheorem{definition}[theorem]{Definition}
\newtheorem{example}[theorem]{Example}

\newtheorem{proposition}[theorem]{Proposition}

\newenvironment{proof}[1][Proof]{\noindent\textbf{#1.} }{\ \rule{0.5em}{0.5em}}
\begin{document}

\title{Spatial pseudoanalytic functions arising from the factorization of linear
second order elliptic operators }
\author{Vladislav V. Kravchenko$^{1}$, S\'{e}bastien Tremblay$^{2}$\\$^{1}$Departamento de Matem\'{a}tica, CINVESTAV del IPN, Unidad\\Quer\'{e}taro, Libramiento Norponiente No.~2000 C.P. 76230 Fracc.\\Real de Juriquilla, Quer\'{e}taro, Mexico\\$^{2}$D\'{e}partement de math\'{e}matiques et d'informatique, Universit\'{e} du\\Qu\'{e}bec, Trois-Rivi\`{e}res, Qu\'{e}bec, G9A 5H7, Canada}
\maketitle

\begin{abstract}
Biquaternionic Vekua-type equations arising from the factorization of linear
second order elliptic operators are studied. Some concepts from classical
pseudoanalytic function theory are generalized onto the considered spatial
case. The derivative and antiderivative of a spatial pseudoanalytic function
are introduced and their applications to the second order elliptic equations
are considered.

\end{abstract}

\section{Introduction}

Theory of pseudoanalytic functions \cite{Berskniga}, \cite{Vekua} is one of
the classical branches of complex analysis extending the concepts and ideas from
analytic function theory onto a much more general situation and involving
linear elliptic equations and systems with variables coefficients. Its
development in forties-fifties of the last century was fast and deep and
historically represented an important impulse to the progress in the general
theory of elliptic systems. Nevertheless the important obstacles for a further
development of pseudoanalytic function theory were its limited practical
applications together with the fact that many important results remained in
the level of existence without a possibility to make them really applicable
for solving problems of mathematical physics. Recent progress in
pseudoanalytic function theory reported in \cite{APFT} shows its deep relation
to the stationary Schr\"{o}dinger equation and more general linear second
order elliptic equations and includes new results which allow one to make the
basic objects of the theory and their applications fully explicit. Among other
results it is worth mentioning the possibility to obtain complete systems of
solutions to second order elliptic equations with variable coefficients and
their use for solving related boundary and eigenvalue problems \cite{CKR}.

The generalization of the concept of a pseudoanalytic function onto the case
of three or more dimensions was a subject of considerable efforts (see, e.g., \cite{APFT},
\cite{Berglez}, \cite{Mal98}, \cite{KR2010}, \cite{Sprossig2001}%
). The development of a spatial generalization of the theory of pseudoanalytic
functions seems to be possible if only the study is restricted to linear first
order systems of relatively special forms otherwise it is difficult to expect
an efficient generalization of the very first concepts of the theory like
derivative and antiderivative not to mention formal powers, the Cauchy
integral formula and other notions. In the present work we study the Vekua
equation arising from the factorization of the stationary Schr\"{o}dinger
operator. This biquaternionic Vekua equation was first introduced in
\cite{K2006}, \cite{APFT} as a direct generalization of the so-called main
Vekua equation closely related to the two-dimensional stationary
Schr\"{o}dinger equation. We show that it is possible in a natural way to
introduce the concept of a derivative for solutions of the biquaternionic main
Vekua equation and exactly as in the two-dimensional case the derivative is a
solution of another Vekua equation the operator of which is involved in the
factorization of the stationary Schr\"{o}dinger operator. Moreover, we also
show that this operation of differentiation is invertible and introduce the
notion of the antiderivative. Among other applications these constructions
allow us to construct exact solutions to the stationary Schr\"{o}dinger
equation from the known ones possessing certain symmetries.

In this work we consider the elliptic case only, however the generalizations
of the presented results onto the hyperbolic situation (see \cite{KRT} and
\cite{KR2010}) are also possible.

\section{Some definitions and results from two-dimensional pseudoanalytic
function theory}

In this section we give some preliminary definitions and results from the
classical (complex) pseudoanalytic function theory. For more details, see
\cite{Berskniga, APFT, BersStat}.

Let $\Omega$ be a simply connected domain in
$\mathbb{R}^{2}$. We say that a pair of complex functions $F$ and $G$
possessing in $\Omega$ partial derivatives with respect to the real variables
$x$ and $y$ is a generating pair if it satisfies the inequality
\[
\operatorname{Im}(\overline{F}G)>0\qquad\text{in }\Omega.
\]
From this definition it follows that if $W:\Omega\rightarrow\mathbb{C}$, then
there exist unique real valued functions $\phi,\psi$ such that
\[
W(z)=\phi(z)F(z)+\psi(z)G(z),\ \ \ \forall z\in\Omega,
\]
where $\phi=\displaystyle \frac{\mathrm{Im}(\overline{W}G)}{\mathrm{Im}%
(\overline{F}G)}$ and $\psi=-\displaystyle \frac{\mathrm{Im}(\overline{W}%
F)}{\mathrm{Im}(\overline{F}G)}$.

For a given generating pair $(F,G)$ in $\Omega$, we can define an
$(F,G)$-derivative $\dot{W}(z_{0})$ on $W:\Omega\rightarrow\mathbb{C}$ if the
(finite) limit
\[
\dot{W}(z_{0})=\lim_{z\rightarrow z_{0}}\frac{{W(z)-\phi(z_{0})F(z)-\psi
(z_{0})G(z)}}{{z-z_{0}}}%
\]
exists. Sometimes instead of the notation $\dot{W}$ for the $(F,G)$-derivative
of $W$ we will use the notation $\displaystyle\frac{d_{(F,G)}W}{dz}$.

We denote by $\partial_{\overline{z}}=\frac{1}{2}(\partial_{x}+i\partial_{y})$
and $\partial_{z}=\frac{1}{2}(\partial_{x}-i\partial_{y})$, where
$\partial_{x}:=\frac{\partial}{\partial x}$ and $\partial_{y}:=\frac{\partial
}{\partial y}$. The following expressions are known as characteristic
coefficients of the pair $(F,G)$
\[
a_{(F,G)}=-\frac{\overline{F}G_{\overline{z}}-F_{\overline{z}}\overline{G}%
}{F\overline{G}-\overline{F}G},\qquad b_{(F,G)}=\frac{FG_{\overline{z}%
}-F_{\overline{z}}G}{F\overline{G}-\overline{F}G},
\]

\[
A_{(F,G)}=-\frac{\overline{F}G_{z}-F_{z}\overline{G}}{F\overline{G}%
-\overline{F}G},\qquad B_{(F,G)}=\frac{FG_{z}-F_{z}G}{F\overline{G}%
-\overline{F}G},
\]
where the subindex $\overline{z}$ or $z$ means the application of
$\partial_{\overline{z}}$ or $\partial_{z}$ respectively.

\begin{theorem}
Let $(F,G)$ be a generating pair in $\Omega$ and $W\in C^{1}(\Omega)$. The
$(F,G)$-derivative $\dot{W}=\displaystyle\frac{d_{(F,G)}W}{dz}$ of $W$ exists
and has the form
\begin{equation}
\dot{W}=\phi_{z}F+\psi_{z}G=W_{z}-A_{(F,G)}W-B_{(F,G)}\overline{W}
\label{derivative}%
\end{equation}
if and only if
\begin{equation}
W_{\overline{z}}=a_{(F,G)}W+b_{(F,G)}\overline{W}. \label{VekuaGen}%
\end{equation}

\end{theorem}

Equation (\ref{VekuaGen}) is called Vekua equation (sometimes Carleman-Vekua
equation) and can also be written in the following form
\[
\phi_{\overline{z}}F+\psi_{\overline{z}}G=0.
\]
The Vekua equation (\ref{VekuaGen}) represents a generalization of the
Cauchy-Riemann system and solutions of this equation are called $(F,G)$%
-pseudoanalytic functions. It is easy to see that functions $F$ and $G$ are
$(F,G)$-pseudoanalytic, and $\overset{\cdot}{F}\equiv\overset{\cdot}{G}%
\equiv0$.

Let $(F,G)$ and $(F_{1},G_{1})$ be two generating pairs in $\Omega$. We say
that $(F_{1},G_{1})$ is a successor of $(F,G)$ and $(F,G)$ is a predecessor of
$(F_{1},G_{1})$ if%
\[
a_{(F_{1},G_{1})}=a_{(F,G)}\qquad\text{and}\qquad b_{(F_{1},G_{1})}%
=-B_{(F,G)}\text{.}%
\]
The importance of this definition becomes obvious from the following result.

\begin{theorem}
\label{ThBersDer}Let $W$ be an $(F,G)$-pseudoanalytic function and let
$(F_{1},G_{1})$ be a successor of $(F,G)$. Then $\overset{\cdot}{W}$ is an
$(F_{1},G_{1})$-pseudoanalytic function.
\end{theorem}

In other words, when $W$ is an $(F,G)$-pseudoanalytic function then $\dot{W}$
satisfies the Vekua equation $(\dot{W})_{\overline{z}}=a_{(F,G)}\dot
{W}-B_{(F,G)}\overline{\dot{W}}$.

For a given generating pair $(F,G)$ in $\Omega$, the adjoint generating pair
$(F,G)^{\ast}=(F^{\ast},G^{\ast})$ is defined by the formulas%
\[
F^{\ast}=-\frac{2\overline{F}}{F\overline{G}-\overline{F}G}\qquad\text{ and }
\qquad G^{\ast}=\frac{2\overline{G}}{F\overline{G}-\overline{F}G}.
\]

The $(F,G)$-integral is defined as follows
\[
\int_{\Gamma}Wd_{(F,G)}z= F(z_{1})\operatorname{Re}\int_{\Gamma}G^{\ast
}Wdz+G(z_{1})\operatorname{Re}\int_{\Gamma}F^{\ast}Wdz
\]
where $\Gamma$ is a rectifiable curve leading from $z_{0}$ to $z_{1}$.

If $W=\phi F+\psi G$ is an $(F,G)$-pseudoanalytic function where $\phi$ and
$\psi$ are real valued functions then
\[
\int_{z_{0}}^{z}\overset{\cdot}{W}d_{(F,G)}z=W(z)-\phi(z_{0})F(z)-\psi
(z_{0})G(z), \label{FGAnt}%
\]
and as $\overset{\cdot}{F}=\overset{}{\overset{\cdot}{G}=}0$, this integral is
path-independent and represents the $(F,G)$-antiderivative of $\overset{\cdot
}{W}$.

Consider the equation
\[
\varphi_{z}=\Phi
\]
for a real-valued function $\varphi$ in a simply connected domain, where $\Phi=\Phi_{1}+i\Phi
_{2}$ is a given complex valued function such that its real part $\Phi_{1}$
and imaginary part $\Phi_{2}$ satisfy the equation
\begin{equation}
\partial_{y}\Phi_{1}+\partial_{x}\Phi_{2}=0, \label{casirot}%
\end{equation}
then we can reconstruct $\varphi$ up to an arbitrary real constant $c$ in the
following way%
\begin{equation}
\varphi(x,y)=2\left(  \int_{\Gamma}\Phi_{1}dx-\Phi_{2}dy\right)  +c.
\label{Antigr2}%
\end{equation}

By $A$ we denote the integral operator in (\ref{Antigr2}):%
\[
A[\Phi](x,y)=2\left(  \int_{\Gamma}\Phi_{1}dx-\Phi_{2}dy\right)  +c.
\]

Thus if $\Phi$ satisfies (\ref{casirot}), there exists a family of real valued
functions $\varphi$ such that $\varphi_{z}=\Phi$, given by the formula
$\varphi=A[\Phi]$.

In a similar way we define the operator $\overline{A}$ corresponding to
equation $\varphi_{\overline{z}}=\Phi$, where $\Phi=\Phi_{1}+i\Phi_{2}$
satisfies $\partial_{y}\Phi_{1}-\partial_{x}\Phi_{2}=0$, by the following
antiderivative operator:
\[
\overline{A}[\Phi](x,y)=2\left(  \int_{\Gamma}\Phi_{1}dx+\Phi_{2}dy\right)
+c.
\]

Now, it is well known that if $f_{0}$ is a non vanishing particular solution of the
one-dimensional stationary Schr\"{o}dinger equation $\big(-\frac{d^{2}}%
{dx^{2}}+q(x)\big)f(x)=0$, the Schr\"{o}dinger operator can be factorized as
$\frac{d^{2}}{dx^{2}}-q=(\frac{d}{dx}+\frac{f_{0}^{\prime}}{f_{0}})(\frac
{d}{dx}-\frac{f_{0}^{\prime}}{f_{0}})$. A similar factorization can be
obtained in the two-dimensional case. First, we define the Vekua operators
\[
V:=\partial_{\overline{z}}-\frac{f_{\overline{z}}}{f}C\quad\text{and}\quad
V_{1}:=\partial_{\overline{z}}+\frac{f_{z}}{f}C,
\]
where $f\in C^{2}(\Omega;\mathbb{R})$ is a nonvanishing real-valued function
defined in $\Omega$ and $C$ is the complex conjugate operator. The
two-dimensional stationary Schr\"{o}dinger equation
\begin{equation}
\big(-\triangle+q(x,y)\big)\varphi=0\qquad\text{ in }\Omega, \label{Schro1}%
\end{equation}
where $\triangle$ is the two-dimensional Laplacian $\triangle:=\frac
{\partial^{2}}{\partial x^{2}}+\frac{\partial^{2}}{\partial y^{2}}$, $q$ and
$\varphi$ are real-valued functions, $\varphi\in C^{2}(\Omega;\mathbb{R})$,
can be factorized as follows.

\begin{theorem}
\emph{\cite{K2005}} Let $f$ be a real-valued nonvanishing solution in $\Omega$
of the Schr\"{o}dinger equation (\ref{Schro1}). Then for any real-valued
function $\varphi\in C^{2}(\Omega)$ the following equalities hold
\begin{equation}
\frac{1}{4}(\triangle-q)\varphi= V_{1}\overline{V}\varphi=\overline{V}%
_{1}V\varphi. \label{factorization}%
\end{equation}

\end{theorem}

\begin{theorem}
\emph{\cite{K2006}} Let $W=W_{1}+iW_{2}$ be a solution of the Vekua equation
$VW=0$. Then $W_{1}=\text{\emph{Re} }W$ is a solution of the Schr\"{o}dinger
equation (\ref{Schro1}) in $\Omega$ and $W_{2}=\text{\emph{Im} }W$ is a
solution of the associated Schr\"{o}dinger equation
\begin{equation}
(-\triangle+r)\psi=0\quad\text{in }\Omega,\quad\quad r=2\left(  \frac{\nabla
f}{f}\right)  ^{2}-q, \label{Schro2}%
\end{equation}
where $(\nabla f)^{2}=f_{x}^{2}+f_{y}^{2}$. \label{VekuasolSchr}
\end{theorem}

Theorem \ref{VekuasolSchr} tells us that as much as real and imaginary parts
of a complex analytic function are harmonic functions, the real and imaginary
parts of a solution of the Vekua equation $VW=0$ are solutions of associated
stationary Schr\"{o}dinger equations. Now, we know that given an arbitrary
real-valued harmonic function in a simply connected domain, a conjugate
harmonic function can be constructed explicitly such that the obtained pair
of harmonic functions represent the real and imaginary parts of a complex
analytic function. This corresponds to the more general situation for
solutions of associated stationary Schr\"{o}dinger equations; the following
theorem gives us the precise result.

\begin{theorem}
\emph{\cite{K2005}\label{ThConjHarm}} Let $W_{1}$ be a real-valued solution of
(\ref{Schro1}) in a simply connected domain $\Omega$. Then the real-valued
function $W_{2}$, a solution of (\ref{Schro2}) such that $W=W_{1}+iW_{2}$ is a
solution of $VW=0$, is constructed according to the formula
\[
W_{2}=f^{-1}\overline{A}[if^{2}\partial_{\overline{z}}(f^{-1}W_{1})].
\]

Given a solution $W_{2}$ of (\ref{Schro2}), the corresponding solution $W_{1}$
of (\ref{Schro1}) such that $W=W_{1}+iW_{2}$ is a solution of $VW=0$, is
constructed as
\[
W_{1}=-f\overline{A}[if^{-2}\partial_{\overline{z}}(fW_{2})].
\]

\end{theorem}

\section{Pseudoanalytic derivative for a class of Vekua equations in the
plane\label{SectTwoDim}}

For $f\in C^{2}(\Omega,\mathbb{R})$, a nonvanishing real-valued function
defined in $\Omega$, it is known \cite[Sect. 3.4]{APFT} that the functions
$F=f$ and $G=i/f$ represent a generating pair for the Vekua equation
\begin{equation}
VW=0, \label{VekuaMain}%
\end{equation}
and if one fixes this generating pair $(F,G)=(f,i/f)$, the equation
\begin{equation}
V_{1}w=0 \label{Vekua1}%
\end{equation}
is satisfied by the $(F,G)$-derivatives of solutions of (\ref{VekuaMain}) and
in particular by the successor $(F_{1},G_{1})$, and hence it is natural to
call $V_{1}$ a successor of $V$.

Note that from (\ref{derivative}) the $(F,G)$-derivative in this case is given
as
\[
\overset{\cdot}{W}=\frac{d_{(F,G)}W}{dz}=W_{z}-\frac{f_{z}}{f}\overline
{W}=CVCW.
\]
Thus, the fact that $\overset{\cdot}{W}$ is a solution of (\ref{Vekua1}) can
be written in the following form. If $W$ is a solution of (\ref{VekuaMain})
then
\[
V_{1}CVCW=0.
\]
Note the following commutation relation%
\[
CV=\overline{V}C
\]
where $\overline{V}=\partial_{z}-\frac{f_{z}}{f}C$. That is,
\[
\overset{\cdot}{W}=\overline{V}W
\]
and $V_{1}\overset{\cdot}{W}=V_{1}\overline{V}W$. The equality to zero of this
last expression in the case when $W$ is a solution of (\ref{VekuaMain}) is a
direct consequence of the factorization (\ref{factorization}), i.e.
$V_{1}\overline{V}\varphi=\frac{1}{4}(\Delta-q)\varphi$ where $q=\frac{\Delta
f}{f}$ and $\varphi$ is real-valued, twice continuously differentiable
function, and
\[
V_{1}\overline{V}(i\phi)=\frac{i}{4}(\Delta-r)\phi
\]
where $r=-q+2\left(  \frac{\nabla f}{f}\right)  ^{2},$ as well as of the fact
that from theorem \ref{VekuasolSchr} if $W$ is a solution of (\ref{VekuaMain})
then $(-\Delta+q)W_{1}=0$ and $(-\Delta+r)W_{2}=0$ where $W_{1}%
=\operatorname{Re}W$ and $W_{2}=\operatorname{Im}W$.

Let us notice that the $(F,G)$-derivative of $W$ can be also represented as
follows%
\begin{equation}
\overline{V}W=f\partial_{z}(f^{-1}W_{1})+if^{-1}\partial_{z}(fW_{2}).
\label{SecondRep1}%
\end{equation}
This form suggests the definition of the inverse operation, the
antiderivative:%
\[
w=fA(f^{-1}\Phi_{1})+if^{-1}A(f\Phi_{2})
\]
where $\Phi=\Phi_{1}+i\Phi_{2}$ is any solution of (\ref{Vekua1}). The
resulting function $w$ is necessarily a solution of (\ref{VekuaMain}).

From (\ref{SecondRep1}) it can be seen that the $(F,G)$-derivative can be
written as
\begin{equation}
\overset{\cdot}{W}=\varphi_{z}F+\psi_{z}G \label{SecondRep}%
\end{equation}
where $\varphi=f^{-1}W_{1}$ and $\psi=fW_{2}$.

\section{Three-dimensional pseudoanalytic function\newline theory for a class
of Vekua equations}

\subsection{General definitions and results}

We denote by $\mathbb{H}(\mathbb{C})=\{Q\ |\ Q=\sum_{\alpha=0}^{3}Q_{\alpha
}\mathbf{e}_{\alpha}\}$ the algebra of complex quaternions (biquaternions)
(see, e.g., \cite{AQA, IRSM}), where $Q_{\alpha}\in\mathbb{C}$, $\mathbf{e}%
_{0}=1$, $\{\mathbf{e}_{k}\ |\ k=1,2,3\}$ are the standard quaternionic
imaginary units satisfying $\mathbf{e}_{\alpha}\mathbf{e}_{\beta}%
+\mathbf{e}_{\beta}\mathbf{e}_{\alpha}=-2\delta_{\alpha,\beta}$, where
$\delta_{\alpha,\beta}$ is the usual Kronecker delta. The imaginary unit in
$\mathbb{C}$ is denoted by $i$ as usual and commutes with all basic units
$\mathbf{e}_{\alpha}$, $\alpha=0,1,2,3$. The vector representation of
$Q\in\mathbb{H}(\mathbb{C})$ can also be used: $Q=\text{Sc}(Q)+\text{Vec}(Q)$,
where $\text{Sc}(Q)=Q_{0}$ and $\text{Vec}(Q)=\mathbf{Q}=\sum_{k=1}^{3}%
Q_{k}\mathbf{e}_{k}$. The quaternionic conjugation is defined by $\overline
{Q}=Q_{0}-\mathbf{Q}$ (sometimes the quaternionic conjugation operator $C_{H}$
is used, i.e. $C_{H}Q=\overline{Q}$).

By $M^{P}$ we denote the operator of multiplication by a biquaternion $P$ from
the right-hand side, i.e.
\[
M^{P}Q=Q\cdot P.
\]

Let $Q$ be a complex quaternion-valued differentiable function of
$\mathbf{x}=(x_{1},x_{2},x_{3})$. We define the Dirac operator (sometimes
called the Moisil-Theodorescu operator) as
\begin{equation}
DQ=\sum_{k=1}^{3} \mathbf{e}_{k}\ \partial_{k} Q, \label{Diracoperator}%
\end{equation}
where $\partial_{k}:=\displaystyle \frac{\partial}{\partial x_{k}}$.
Expression (\ref{Diracoperator}) can be rewritten in a vector form as
\[
DQ=-\text{div } \mathbf{Q}+\nabla Q_{0}+\text{rot } \mathbf{Q},
\]
in other words $\text{Sc}(DQ)=-\text{div }\mathbf{Q}$ and $\text{Vec}%
(DQ)=\nabla Q_{0} +\text{rot } \mathbf{Q}$. Let us notice that $D^{2}%
=-\triangle$, where $\triangle:=\sum_{k=1}^{3} \partial_{k}^{2}$ is the
Laplacian in $\mathbb{R}^{3}$.

We notice that the Dirac operator was introduced as acting from the left-hand
side. The corresponding operator acting from the right-hand side we will
denote by $D_{r}Q=\sum_{k=1}^{3}\partial_{k}Q\ \mathbf{e}_{k}$. In the vector
form $D_{r}$ can be represented as $D_{r}Q=-\text{div }\mathbf{Q}+\nabla
Q_{0}-\text{rot }\mathbf{Q}$.

The Dirac operator admits the following generalization of the Leibniz rule
(see, e.g., \cite{GurlebeckSprossig}):

\begin{theorem}
Let $\{P,Q\}\subset C^{1}\big(\Omega;\mathbb{H}(\mathbb{C})\big)$, where
$\Omega$ is some domain in $\mathbb{R}^{3}$. Then
\[
D[P\cdot Q]=D[P]\cdot Q+\overline{P}\cdot D[Q]+2\big(\text{\emph{Sc}%
}(PD)\big)[Q],
\]
where
\[
\big(\text{\emph{Sc}}(PD)\big)[Q]=-\sum_{k=1}^{3} P_{k}\partial_{k} Q.
\]

\end{theorem}

In particular, we observe that when $\text{Vec}(P)=0$, i.e. $P=P_{0}$, then
\[
D[P_{0}\cdot Q]=D[P_{0}]\cdot Q+P_{0}\cdot D[Q].
\]

Consider now the equation
\[
\nabla\varphi=\mathbf{\Psi}
\]
where $\mathbf{\Psi}=\sum_{k=1}^{3} \Psi_{k} \mathbf{e}_{k}$ is a given
complex-valued vector such that $\text{rot }\mathbf{\Psi}\equiv0$. The
complex-valued scalar function $\varphi$ is then said to be the potential of
$\mathbf{\Psi}$. We will write $\varphi=\mathcal{A}[\mathbf{\Psi}]$. The
operator $\mathcal{A}$ is a simple generalization of the opeators $A$ and
$\overline{A}$ of the preceding section. It is well known that we can
reconstruct $\varphi$, up to an arbitrary complex constant $c$, in the
following way
\[
\mathcal{A}[\mathbf{\Psi}](x,y,z)=\int_{x_{0}}^{x}\Psi_{1}(\xi,y_{0}%
,z_{0})d\xi+ \int_{y_{0}}^{y}\Psi_{2}(x,\eta,z_{0})d\eta+\int_{z_{0}}^{z}%
\Psi_{3}(x,y,\zeta)d\zeta+c,
\]
where $(x_{0},y_{0},z_{0})$ is an arbitrary point in the domain of interest.

Let us now consider the three-dimensional stationary Schr\"{o}dinger equation
\begin{equation}
(-\Delta+q)\phi=0\qquad\text{ in }\Omega, \label{Schro3d}%
\end{equation}
where $q$ and $\phi$ are complex-valued functions depending on $\mathbf{x}%
=(x_{1},x_{2},x_{3})$ and $\Omega$ is a domain in $\mathbb{R}^{3}$. Using the
Dirac operator, one can factorize this equation in a way similar to that in
theorem~\ref{factorization}.

\begin{theorem}
\emph{\cite[Sect. 16.1]{APFT}} Let $f$ be a nonvanishing particular solution
of (\ref{Schro3d}). Then for any scalar (complex-valued) function $\phi\in
C^{2}(\Omega;\mathbb{C})$ the following equality holds:
\begin{equation}
(-\Delta+q)\phi=(D+M^{\frac{Df}{f}})(D-\frac{Df}{f}C_{H})\phi. \label{Fact}%
\end{equation}

\end{theorem}

The equivalent of theorem \ref{VekuasolSchr} in the three-dimensional case is
given here.

\begin{theorem}
\label{ThSolutionsOfVekuaMain3D}Let $W_{0}+\mathbf{W}$ be a solution of the
equation
\begin{equation}
\left(  D-\frac{Df}{f}C_{H}\right)  W=0. \label{mainVekuaquat}%
\end{equation}
Then $W_{0}$ is a solution of (\ref{Schro3d}) with $q=\frac{\triangle f}{f}$;
the scalar function $\varphi_{0}=W_{0}/f$ is a solution of the equation
\begin{equation}
\text{\emph{div}}(f^{2}\nabla\varphi_{0})=0 \label{eqphi0}%
\end{equation}
and the vector function $\mathbf{\Phi}=f\mathbf{W}$ is a solution of the
equation
\begin{equation}
\text{\emph{rot}}(f^{-2}\text{ \emph{rot} }\mathbf{\Phi})=0. \label{eqPhi}%
\end{equation}

\end{theorem}

\subsection{Pseudoanalytic derivative in the space for a class of Vekua
equations}

We denote the operators
\[
V:=D-\frac{Df}{f}C_{H}\quad\text{and}\quad\overline{V}_{1}:=D+M^{\frac{Df}{f}%
}\text{.}%
\]
Here the notations from Section \ref{SectTwoDim} are preserved in order to
make the analogy between the two-dimensional and the higher dimensional case
more transparent.

Application of the operator $C_{H}$ to (\ref{Fact}) gives us
\[
(-\Delta+q)\phi=(D_{r}+\frac{Df}{f})(D_{r}-M^{\frac{Df}{f}}C_{H})\phi.
\]
We denote
\[
\overline{V}:=D_{r}-M^{\frac{Df}{f}}C_{H}\quad\text{and}\quad V_{1}%
:=D_{r}+\frac{Df}{f}.
\]
All the operators $V$, $\overline{V}$, $V_{1}$ and $\overline{V}_{1}$ are
considered in application to $\mathbb{H}(\mathbb{C})$-valued, continuously
differentiable functions. Note that the following relations hold%
\[
C_{H}V=-\overline{V}C_{H}\quad\text{and}\quad C_{H}V_{1}=-\overline{V}%
_{1}C_{H}.
\]

Suppose
\begin{equation}
VW=0 \label{VekuaMain3D}%
\end{equation}
where $W\in C^{1}(\Omega;\mathbb{H}(\mathbb{C}))$. Then due to (\ref{Fact}) we
have
\begin{equation}
\overline{V}_{1}VW_{0}=V_{1}\overline{V}W_{0}=(-\Delta+q)W_{0}=0
\label{factoriz}%
\end{equation}
and hence%
\[
0=\overline{V}_{1}V\mathbf{W}=C_{H}\overline{V}_{1}V\mathbf{W}=V_{1}%
\overline{V}\mathbf{W}.
\]
Thus, from this and from (\ref{factoriz}) the $\mathbb{H}(\mathbb{C})$-valued
function $\overline{V}W$ is a solution of the equation
\begin{equation}
V_{1}w=0 \label{Vekua1_3D}%
\end{equation}
whenever $W$ is a solution of (\ref{VekuaMain3D}). Consequently we define the
Bers derivative of a solution $W$ of (\ref{VekuaMain3D}) as follows

\begin{definition}
Let $W$ be a solution of (\ref{VekuaMain3D}). Then the function
\[
\overset{\cdot}{W}=\overline{V}W
\]
is a solution of (\ref{Vekua1_3D}) and is called the derivative of $W$.
\end{definition}

In other words, $\overline{V}:\ker V\rightarrow\ker V_{1}$. Moreover, as we
show below the derivative $\overset{\cdot}{W}$ is purely vectorial, and hence
it is a solution of the equation $\overline{V}_{1}w=0$ as well.

Let us obtain another representation for $\overset{\cdot}{W},$ similar to
(\ref{SecondRep}). We have%
\[
\overline{V}W=\left(  D_{r}-M^{\frac{Df}{f}}C_{H}\right)  W=\left(
D_{r}-M^{\frac{Df}{f}}C_{H}\right)  \sum_{\alpha=0}^{3}\varphi_{\alpha
}F_{\alpha}%
\]
where $F_{\alpha}$ represent the generating quartet \cite[Sect. 16.2]{APFT}
for equation (\ref{VekuaMain3D})%
\[
F_{0}=f,\quad F_{1}=\frac{\mathbf{e}_{1}}{f},\quad F_{2}=\frac{\mathbf{e}_{2}%
}{f},\quad F_{3}=\frac{\mathbf{e}_{3}}{f},
\]
$\varphi_{k}$ are scalar functions defined by $\varphi_{0}=W_{0}/f$ and
$\varphi_{k}=fW_{k}$ for $k=1,2,3$. Then%
\begin{align*}
\overline{V}W  &  =\sum_{\alpha=0}^{3}F_{\alpha}\cdot D\varphi_{\alpha}%
+\sum_{\alpha=0}^{3}\varphi_{\alpha}\cdot D_{r}F_{\alpha}-\sum_{\alpha=0}%
^{3}\varphi_{\alpha}\overline{F}_{\alpha}\frac{Df}{f}\\
&  =\sum_{\alpha=0}^{3}F_{\alpha}\cdot D\varphi_{\alpha}+\sum_{\alpha=0}%
^{3}\varphi_{\alpha}\cdot\overline{V}F_{\alpha}.
\end{align*}
It is easy to see that $\overline{V}F_{\alpha}=0$ for $\alpha=0,1,2,3$.
Indeed, we have
\[
0=VF_{\alpha}=C_{H}VF_{\alpha}=\overline{V}\overline{F}_{\alpha}%
\]
and if $\alpha=0$, then $\overline{V}\overline{F}_{0}=\overline{V}F_{0}=0$,
otherwise $\overline{V}\overline{F}_{k}=-\overline{V}F_{k}=0$ for $k=1,2,3$.
Thus, another representation for the derivative $\overset{\cdot}{W}$ has the
form
\begin{equation}
\overset{\cdot}{W}=\overline{V}W=\sum_{\alpha=0}^{3}F_{\alpha}\cdot
D\varphi_{\alpha}. \label{Derivative}%
\end{equation}
Let us obtain the inverse operation,- the antiderivative. We have the
following two equalities%
\begin{equation}
\sum_{\alpha=0}^{3}\left(  D\varphi_{\alpha}\right)  F_{\alpha}=0
\label{VekuaSecond}%
\end{equation}
and
\begin{equation}
\sum_{\alpha=0}^{3}\left(  D\varphi_{\alpha}\right)  \overline{F}_{\alpha
}=-\overline{\overset{\cdot}{W}}, \label{VekuaDer}%
\end{equation}
the first of them being another form of the Vekua equation (\ref{VekuaMain3D})
(see \cite[Sect. 16.2]{APFT}), and the second is obtained from
(\ref{Derivative}) by applying $C_{H}$. More explicitly,%
\[
D\varphi_{0}\cdot f+\sum_{k=1}^{3}D\varphi_{k}\cdot\frac{\mathbf{e}_{k}}{f}=0
\]
and
\[
D\varphi_{0}\cdot f-\sum_{k=1}^{3}D\varphi_{k}\cdot\frac{\mathbf{e}_{k}}%
{f}=-\overline{\overset{\cdot}{W}}.
\]
Then
\begin{equation}
D\varphi_{0}=-\frac{1}{2f}\overline{\overset{\cdot}{W}}. \label{Dphi0}%
\end{equation}
A scalar solution $\varphi_{0}$ of this equation exists if only $\overset
{\cdot}{W}$ is purely vectorial and satisfies the equation
\[
\operatorname*{rot}\left(  \frac{1}{f}\overset{\cdot}{W}\right)  =0.
\]
Consider
\begin{align*}
\overset{\cdot}{W}+\overline{\overset{\cdot}{W}}  &  =\sum_{\alpha=0}%
^{3}F_{\alpha}\cdot D\varphi_{\alpha}-\sum_{\alpha=0}^{3}\left(
D\varphi_{\alpha}\right)  \overline{F}_{\alpha}\\
&  =\sum_{k=1}^{3}\frac{\mathbf{e}_{k}}{f}D\varphi_{k}+\sum_{k=1}^{3}%
D\varphi_{k}\cdot\frac{\mathbf{e}_{k}}{f}\\
&  =\frac{1}{f}\sum_{k=1}^{3}\left(  \mathbf{e}_{k}D\varphi_{k}+D\varphi
_{k}\mathbf{e}_{k}\right) \\
&  =-\frac{1}{f}\sum_{k=1}^{3}\partial_{k}\varphi_{k}=-\frac{1}{f}%
\operatorname*{div}\mathbf{\Phi}%
\end{align*}
where $\mathbf{\Phi}:=\sum_{k=1}^{3}\varphi_{k}\mathbf{e}_{k}$.

As $\mathbf{\Phi}=f\mathbf{W}$ we have that \cite[equation (16.11)]{APFT}
\begin{equation}
\operatorname*{div}\mathbf{\Phi}=0. \label{divPhi}%
\end{equation}
Thus, $\overset{\cdot}{W}+\overline{\overset{\cdot}{W}}=0$ and hence
$\overset{\cdot}{W}$ is purely vectorial.

Consider the expression
\begin{align*}
\operatorname*{rot}\left(  \frac{1}{f}\overset{\cdot}{W}\right)   &
=\operatorname*{rot}\left(  \frac{1}{f}\sum_{\alpha=0}^{3}F_{\alpha}\cdot
D\varphi_{\alpha}\right) \\
&  =\operatorname*{rot}\left(  \nabla\varphi_{0}+\frac{1}{f^{2}}\sum_{k=1}%
^{3}\mathbf{e}_{k}D\varphi_{k}\right) \\
&  =\operatorname*{rot}\left(  \frac{1}{f^{2}}\sum_{k=1}^{3}\mathbf{e}%
_{k}D\varphi_{k}\right)  .
\end{align*}
From the previous chain of equalities we already know that
\[
\sum_{k=1}^{3}\mathbf{e}_{k}D\varphi_{k}=-\sum_{k=1}^{3}D\varphi_{k}%
\mathbf{e}_{k}=-D\mathbf{\Phi}=-\operatorname*{rot}\mathbf{\Phi},
\]
because of (\ref{divPhi}). Thus,
\[
\operatorname*{rot}\left(  \frac{1}{f}\overset{\cdot}{W}\right)
=-\operatorname*{rot}\left(  \frac{1}{f^{2}}\operatorname*{rot}\mathbf{\Phi
}\right)
\]
and according to theorem \ref{ThSolutionsOfVekuaMain3D} this is zero (see
\cite[theorem 161]{APFT}).

We have then that if $W$ is a solution of (\ref{VekuaMain3D}) then equation
(\ref{Dphi0}) possesses a solution $\varphi_{0}$ which can be written as
follows%
\begin{equation}
\varphi_{0}=\frac{1}{2}\mathcal{A}\left[  \frac{1}{f}\overset{\cdot}%
{W}\right]  . \label{phi0}%
\end{equation}
This expression is unique up to an additive constant.

Now let us see, how $\mathbf{\Phi}$ can be recovered from (\ref{VekuaSecond})
and (\ref{VekuaDer}). We have
\[
\sum_{k=1}^{3}D\varphi_{k}\cdot\frac{\mathbf{e}_{k}}{f}=\frac{1}{2}%
\overline{\overset{\cdot}{W}},
\]
that is,
\begin{equation}
D\mathbf{\Phi}=\frac{f}{2}\overline{\overset{\cdot}{W}}. \label{DPhi}%
\end{equation}
More explicitly,%
\[
\operatorname*{div}\mathbf{\Phi}=0\quad\text{and}\quad\operatorname*{rot}%
\mathbf{\Phi=-}\frac{f}{2}\overset{\cdot}{W}.
\]
Notice that $\operatorname*{div} \big(f\overset{\cdot}{W}\big)  =0.$
Indeed, according to (\ref{Dphi0}),%
\[
f\overset{\cdot}{W}=2f^{2}\nabla\varphi_{0}=2f^{2}\nabla\left(  \frac{W_{0}%
}{f}\right)
\]
and due to theorem \ref{ThSolutionsOfVekuaMain3D}, $\operatorname*{div}\left(
f^{2}\nabla\left(  \frac{W_{0}}{f}\right)  \right)  =0$.

Thus, the problem of recovering $\mathbf{\Phi}$ reduces to the well-studied
problem of reconstruction of a vector by its divergence and rotor.

We will need the following notation. For a given vector function $\mathbf{Q}$,
we define $\mathbf{B}[\mathbf{Q}](\mathbf{x})$ by
\[
\mathbf{B}[\mathbf{Q}](\mathbf{x})=\frac{1}{4\pi}\int_{\Omega}\frac
{\mathbf{Q}(\mathbf{y})}{\left\vert \mathbf{x}-\mathbf{y}\right\vert }%
d\Omega,
\]
Note that $\mathbf{B}$ is a right-inverse for the operator
$\operatorname*{rot}\operatorname*{rot}$.

Then we have (see, e.g., \cite[Sect. 5.7]{Korn})%
\begin{equation}
\mathbf{\Phi}=\operatorname*{rot}\left(  \mathbf{B}\left[  -\frac{f}{2}\dot
{W}\right]  \right)  +\nabla h, \label{Phi}%
\end{equation}
where $h$ is any harmonic function defined in $\Omega$.

These results allow us to formulate the following theorem.

\begin{theorem}
Let $\mathbf{w}\in C^{1}\big(\Omega;\mathbb{H}(\mathbb{C})\big)$ be a purely
vectorial solution of the equation
\begin{equation}
\left(  D+M^{\frac{Df}{f}}\right)  \mathbf{w}=0 \label{D+M}%
\end{equation}
where $f$ is a nonvanishing scalar $C^{1}$-function. Then a solution of the
equation
\begin{equation}
DW-\frac{Df}{f}\overline{W}=0\qquad\text{in }\Omega\label{VekuaMain3DExplicit}%
\end{equation}
such that $\mathbf{w}=\overset{\cdot}{W}$ is defined as follows%
\begin{equation}
W=\frac{1}{2}\left(  f\mathcal{A}\left[  \frac{\mathbf{w}}{f}\right]
-\frac{1}{f}\operatorname*{rot}\left(  \mathbf{B}\left[  f\mathbf{w}\right]
\right)  +\frac{\nabla h}{f}\right)  \label{Antider}%
\end{equation}
where $h$ is an arbitrary harmonic function in $\Omega$.
\end{theorem}

\begin{proof}
First, we notice that
\begin{equation}
\operatorname*{rot}\left(  \frac{\mathbf{w}}{f}\right)  =0. \label{rot=0}%
\end{equation}
Indeed, from the vector part of (\ref{D+M}) one has $\operatorname*{rot}%
\mathbf{w}+(\mathbf{w}\times\frac{\nabla f}{f})=0$ which is equivalent to
(\ref{rot=0}). Thus, the vector $\frac{\mathbf{w}}{f}$ is a gradient and the
expression $\mathcal{A}\left[  \frac{\mathbf{w}}{f}\right]  $ is well defined.

Observe that the operator in equation (\ref{VekuaMain3DExplicit}) can be
written in the following form
\[
D-\frac{Df}{f}C_{H}=\frac{1}{2}\left(  fDf^{-1}(I+C_{H})+f^{-1}Df(I-C_{H}%
)\right)
\]
where $I$ is the identity operator. Applying it to (\ref{Antider}) we obtain
\begin{align*}
& \frac{1}{4}\left(  fDf^{-1}(I+C_{H})+f^{-1}Df(I-C_{H})\right)  \left(
f\mathcal{A}\left[  \frac{\mathbf{w}}{f}\right]  -\frac{1}{f}%
\operatorname*{rot}\left(  \mathbf{B}\left[  f\mathbf{w}\right]  \right)
+\frac{\nabla h}{f}\right) \\*[2ex]
 = & \frac{1}{2}\left(  fD\mathcal{A}\left[  \frac{\mathbf{w}}{f}\right]
-f^{-1}D\operatorname*{rot}\left(  \mathbf{B}\left[  f\mathbf{w}\right]
\right)  \right) =\frac{1}{2}\left(  \mathbf{w-w}\right)  =0.
\end{align*}
Next we prove that $\mathbf{w}=\overset{\cdot}{W}$. For this we observe that
\begin{align*}
\overset{\cdot}{W}& =\overline{V}W=\frac{1}{2}\left(  fD_{r}f^{-1}%
(I+C_{H})+f^{-1}D_{r}f(I-C_{H})\right)  W \\*[2ex] = &
\frac{1}{2}\left(  fD\mathcal{A}\left[  \frac{\mathbf{w}}{f}\right]
+f^{-1}D\operatorname*{rot}\left(  \mathbf{B}\left[  f\mathbf{w}\right]
\right)  \right)  =\frac{1}{2}\left(  \mathbf{w+w}\right)  =\mathbf{w}.
\end{align*}
\end{proof}

Another important corollary of the preceding calculations is the following
generalization of theorem \ref{ThConjHarm}.

\begin{theorem}
Let $W_{0}$ be a scalar solution of (\ref{Schro3d}) with $q=\frac{\triangle
f}{f}$ in $\Omega$. Then the vector function $\mathbf{W}$ such that
$\mathbf{\Phi}=f\mathbf{W}$ be a solution of (\ref{eqPhi}) and (\ref{divPhi})
and $W=W_{0}+\mathbf{W}$ be a solution of (\ref{mainVekuaquat}), is
constructed according to the formula
\begin{equation}
\mathbf{W}=-f^{-1}\left\{  \text{\emph{rot}}\left(  \mathbf{B}\left[
f^{2}\ \nabla\left(  f^{-1}W_{0}\right)  \right]  \right)  +\nabla h\right\}
. \label{Wvect}%
\end{equation}

Given a solution $\mathbf{\Phi}$ of (\ref{eqPhi}) and (\ref{divPhi}), for
$\mathbf{W=}\frac{1}{f}\mathbf{\Phi}$ the corresponding solution $W_{0}$ of
(\ref{Schro3d}) such that $W=W_{0}+\mathbf{W}$ be a solution of
(\ref{mainVekuaquat}), is constructed as follows
\begin{equation}
W_{0}=-f\mathcal{A}\left[  f^{-2}\operatorname*{rot}(f\mathbf{W})\right]  .
\label{W0}%
\end{equation}

\end{theorem}

\begin{proof}
Let $W_{0}$ be a scalar solution of (\ref{Schro3d}) with $q=\frac{\triangle
f}{f}$. Then from (\ref{Dphi0}) we have that $\overset{\cdot}{W}%
=2f\nabla\left(  f^{-1}W_{0}\right)  $ and due to (\ref{Phi}), we obtain
\[
\mathbf{\Phi}=-\operatorname*{rot}\left(  \mathbf{B}\left[  f^{2}%
\ \nabla\left(  f^{-1}W_{0}\right)  \right]  \right)  +\nabla h,
\]
and as a consequence (\ref{Wvect}) holds.

Now let $\mathbf{\Phi}$ be a solution of (\ref{eqPhi}) and (\ref{divPhi}).
According to (\ref{DPhi}) we have $\overset{\cdot}{W}=-2f^{-1}D\mathbf{\Phi
=}-2f^{-1}\operatorname*{rot}\left(  f\mathbf{W}\right)  $. Now from
(\ref{phi0}) we obtain (\ref{W0}).
\end{proof}

\section{Examples and applications}

In the present section we analize the possibility to construct solutions for
the equation (\ref{D+M}) under certain restrictions on the function $f$ and
apply them for obtaining solutions of the Schr\"{o}dinger equation
(\ref{Schro3d}).

\begin{proposition}
\label{PropParallel}Let $f$ be a nonvanishing $C^{1}$ scalar function defined
in $\Omega$ and there exist a harmonic function $\rho$ such that
\begin{equation}
\nabla f\times\nabla\rho  =0\text{\quad in }\Omega
.\label{parallel}%
\end{equation}
Then the purely vectorial function
\begin{equation}
\mathbf{F}=\frac{1}{f}D\rho\label{F}%
\end{equation}
is a solution of (\ref{D+M}) and $\mathbf{G}=fD\rho$ is a solution of the
equation
\begin{equation}
\left(  D-M^{\frac{Df}{f}}\right)  \mathbf{G}=0.\label{D-M}%
\end{equation}

\end{proposition}

\begin{proof}
Consider $D\mathbf{F}=-\frac{1}{f^{2}}Df\,D\rho=-\frac{1}{f}D\rho\,\frac
{Df}{f}=-\mathbf{F}\frac{Df}{f}$. In a similar way one can check the validity
of (\ref{D-M}).
\end{proof}

Let us notice that (\ref{parallel}) is fulfilled when $f$ depends on $\rho$
and this is of a considerable importance because harmonic functions are behind
all physically interesting orthogonal coordinate systems. For example,
consider the spherical coordinate $r=\sqrt{x^{2}+y^{2}+z^{2}}$ and let
$f=f(r)$. The function $r$ itself is not harmonic, however the function
$\rho=1/r$ is harmonic everywhere except the origin. Thus, when $f$ is a
function of $r$ and hence a function of $\rho$, and the domain $\Omega$ does
not contain the origin the function (\ref{F}) is a solution of (\ref{D+M}).

\begin{proposition}
\label{PropOrthogonal}Let $f$ be a nonvanishing $C^{1}$ scalar function
defined in $\Omega$ and there exist a harmonic function $\rho$ such that
\[
\left\langle \nabla f,\nabla\rho\right\rangle =0\text{\quad in }\Omega.
\]
Then the purely vectorial function
\[
\mathbf{F}=fD\rho
\]
is a solution of (\ref{D+M}) and $\mathbf{G}=\frac{1}{f}D\rho$ is a solution
of (\ref{D-M}).
\end{proposition}

\begin{proof}
Consider $D\mathbf{F}=Df\,D\rho=-fD\rho\,\frac{Df}{f}=-\mathbf{F}\frac{Df}{f}%
$. In a similar way one can check the validity of (\ref{D-M}).
\end{proof}

These two propositions give us the possibility to obtain three independent
solutions of equation (\ref{D+M}) in a simple way when $f$ is a function of a
harmonic function and this harmonic function possesses two different
orthogonal harmonic functions. As we show in the next example this situation
is of practical importance. Let us notice that in general for a given harmonic
function in three dimensions there not necessarilly exists a nonconstant
orthogonal harmonic function (a counterexample is given in \cite{Elkins}).

\begin{example}
Consider the cylindrical coordinate system $r=\sqrt{x^{2}+y^{2}}$,
$\theta=\arctan(y/x)$ and $z$. Suppose that $f$ depends only on $r$. Then the
required harmonic function $\rho$ can be chosen as $\rho=\log r$, and
obviously the functions $\theta$ and $z$ are harmonic and orthogonal to $\rho
$. Thus $\mathbf{F}_{1}=\frac{1}{f(r)}D\rho=-\frac{1}{f(r)}\left(  \frac
{x}{r^{3}}\mathbf{e}_{1}+\frac{y}{r^{3}}\mathbf{e}_{2}\right)  $,
$\mathbf{F}_{2}=f(r)D\theta=-f(r)\left(  -\frac{y}{r^{2}}\mathbf{e}_{1}%
+\frac{x}{r^{2}}\mathbf{e}_{2}\right)  $ and $\mathbf{F}_{3}=f(r)\mathbf{e}%
_{3}$ are solutions of (\ref{D+M}) independent in the sense that any function
$\mathbf{w}$ defined in a domain not containing the axis $r=0$ and having
purely quaternionic values can be represented in the form $\mathbf{w}%
=\sum_{k=1}^{3}\varphi_{k}\mathbf{F}_{k}$ where $\varphi_{k}$ are scalar
functions. Consequently $\mathbf{F}_{1}$, $\mathbf{F}_{2}$ and $\mathbf{F}%
_{3}$ form a generating triplet of equation (\ref{D+M}) (see more on this
concept in \cite{KR2010}) and $\mathbf{w}$ is a solution of (\ref{D+M}) if and
only if the scalar functions $\varphi_{k}$ satisfy the equation
\begin{equation}
\sum_{k=1}^{3}D\varphi_{k}\,\mathbf{F}_{k}=0\text{\quad in }\Omega
.\label{Vekua3Dseckind}%
\end{equation}
In general we note that if $f$ is a function of a harmonic function depending
on two Cartesian variables then we always can obtain a generating triplet for
equation (\ref{D+M}) and hence reduce the equation to the form
(\ref{Vekua3Dseckind}). Indeed, for such harmonic function we always have two
orthogonal harmonic functions: its conjugate harmonic depending on two
Cartesian variables as well, and the third Cartesian coordinate.
\end{example}

Assume that $f$ is a function of a harmonic function $\rho$. Then according to
proposition \ref{PropParallel}, $\mathbf{F}=\frac{1}{f}D\rho$ is a solution of
(\ref{D+M}). Its antiderivative (\ref{Antider}) gives us a solution of
(\ref{VekuaMain3DExplicit}) and according to theorem
\ref{ThSolutionsOfVekuaMain3D} the scalar part of the obtained solution is a
solution of (\ref{Schro3d}). Thus we obtain that $\psi=f\mathcal{A}\left[
\frac{\mathbf{F}}{f}\right]  =f\mathcal{A}\left[  \frac{D\rho}{f^{2}}\right]
$ is necessarily a solution of (\ref{Schro3d}), and hence the results of this
section allow us to construct new solutions of the Schr\"{o}dinger equation (\ref{Schro3d}) from known
solutions possessing certain symmetries.

\section{Conclusions}

For biquaternionic Vekua type equations arising from the factorization of
linear second order elliptic operators we introduced and studied the concepts
of the derivative and of the antiderivative in the sense of Bers. These
concepts represent new relations between first order elliptic systems and
allow us to study solutions of the second order elliptic equations in relation
to components of the biquaternion valued solutions of biquaternionic Vekua
equations. We have shown that the revealed relations are useful for obtaining
new exact solutions of the Schr\"{o}dinger type equations and at the same time
they imply natural open questions important for the general theory of linear
second order elliptic equations. Namely, we showed how in a quite general and
important for practical applications situation one can obtain a generating
triplet for the biquaternionic Vekua equation which can be considered as a
successor of the main Vekua equation (\ref{VekuaMain3DExplicit}). In the
two-dimensional case this was sufficient to obtain a whole generating sequence
related to the Vekua equation and as a consequence a complete system of
solutions called formal powers (see \cite{APFT}) for the Vekua equation and
for the related second order elliptic equation. At present it is not clear how
to generalize those results onto the spatial situation.

\subsection*{Acknowledgments}
V.K. acknowledges the support of CONACYT of Mexico via the research project 50424.
The research of S.T. is partly supported by grant from NSERC of Canada.

\end{document}